\documentclass[11pt]{article}
\usepackage{mathrsfs}
\usepackage{epic,eepic,epsf,epsfig}
\usepackage{amsfonts,srcltx,mathrsfs}
\usepackage{multirow}
\usepackage{amsmath}
\usepackage{amssymb}
\usepackage{amsbsy}
\usepackage{graphicx}
\usepackage{amsfonts}%,srcltx}
\usepackage{color}
\usepackage{setspace}

\newtheorem{lem}{Lemma}[section]%
\newtheorem{theorem}[lem]{Theorem}%
\newtheorem{exam}[lem]{Example}%
\newtheorem{prop}[lem]{Proposition}%

\def\a{\alpha}

 \def\O{\Omega} \def\G{\Gamma}

\def\di{\bigm|} \def\lg{\langle} \def\rg{\rangle}
\def\nd{\mathrel{\bigm|\kern-.7em/}}

\def\f{\noindent}

\def\PSL{\hbox{\rm PSL}}\def\PSU{\hbox{\rm PSU}}
 \def\PGL{\hbox{\rm PGL}}  \def\Mult{\hbox{\rm Mult}}
\def\PSp{\hbox{\rm PSp}}\def\P\GammaL{\hbox{\rm P\GammaL}}

\def\Aut{\hbox{\rm Aut}}

\def\Inn{\hbox{\rm Inn}}

\def\Cay{\hbox{\rm Cay}}
\def\Cos{\hbox{\rm Cos}}

\newcommand{\qed}{\mbox{\raisebox{0.7ex}{\fbox{}}} \vspace{4truemm}}
\def\mz{{\mathbb Z}}

\begin{document}
\title{On basic graphs of symmetric graphs of valency five}

\footnotetext[1]{ Corresponding author. E-mails:
yangdawei$@$math.pku.edu.cn (D.-W. Yang), yqfeng$@$bjtu.edu.cn (Y.-Q. Feng), jinkwak@bjtu.edu.cn
(J.H. Kwak), julee@ynu.ac.kr (J. Lee)}

\author{Da-Wei Yang$^{\rm a}$, Yan-Quan Feng$^{\rm a,}$\footnotemark, Jin Ho Kwak$^{\rm a}$, Jaeun Lee$^{\rm b}$\\
{\small\em $^{\rm a}$ Department of Mathematics, Beijing Jiaotong
University, Beijing 100044, China} \\
{\small\em $^{\rm b}$ Mathematics, Yeungnam University, Kyongsan
712-749, Korea}\\ \\ \\ 
{Dedicated to the memory of Michel Deza}\\}

\date{}
 \maketitle

\begin{abstract}

A graph $\G$ is {\em symmetric} or {\em arc-transitive} if its automorphism group
$\Aut(\G)$ is transitive on the arc set of the graph, and $\G$ is {\em basic} if
$\Aut(\G)$ has no non-trivial normal subgroup $N$ such that
the quotient graph $\G_N$ has the same valency with $\G$.
In this paper, we classify symmetric basic graphs
of order $2qp^n$ and valency 5, where $q<p$ are two primes
and $n$ is a positive integer. It is shown that such a graph
is isomorphic to a family of Cayley graphs on dihedral groups
of order $2q$ with $5\di (q-1)$, the complete graph $K_6$ of order $6$,
the complete bipartite graph $K_{5,5}$ of order 10,
or one of the nine sporadic coset graphs associated with non-abelian simple groups.
As an application, connected pentavalent symmetric graphs
of order $kp^n$ for some small integers $k$ and $n$ are classified.

\bigskip
\f {\bf Keywords:} Symmetric graph, arc-transitive, normal cover, basic graph.

\end{abstract}

\section{Introduction}

Throughout this paper, all groups and graphs are finite, and all graphs are simple, undirected
and connected. Let $G$ be a permutation group on a set $\O$ and let $\a\in \O$. Denote by $G_{\a}$ the stabilizer of
$\a$ in $G$, that is, the subgroup of $G$ fixing the point $\a$. We say that $G$ is {\it semiregular} on
$\O$ if $G_{\a}=1$ for every $\a\in\O$, and {\it regular} if $G$ is transitive and semiregular.

For a graph $\G$, we denote its vertex set,
edge set and automorphism group by $V(\G)$, $E(\G)$ and $\Aut(\G)$,
respectively. An {\it $s$-arc} in a graph $\G$ is an ordered
$(s+1)$-tuple $(v_0,v_1,\ldots,v_s)$ of vertices of $\G$ such that
$v_{i-1}$ is adjacent to $v_i$ for $1\leq i\leq s$, and $v_{i-1}\neq
v_{i+1}$ for $1\leq i<s$. A $1$-arc is just called an {\it arc}. A graph $\G$ is said to be
{\it $(G,s)$-arc-transitive} if $G$ acts transitively
on the set of $s$-arcs of $\G$, and {\it $(G, s)$-transitive} if $G$ acts transitively on the $s$-arcs but not on the $(s+1)$-arcs of $\G$,
where $G$ is a subgroup of $\Aut(\G)$.
A graph $\G$ is said to be {\it $s$-arc-transitive} or {\it $s$-transitive} if it is
$(\Aut(\G), s)$-arc-transitive or $(\Aut(\G), s)$-transitive, respectively.
In particular, 0-arc-transitive means {\it vertex-transitive}, and 1-arc-transitive means {\it arc-transitive} or {\it symmetric}.

Let $\G$ be a graph and $N\leq \Aut(\G)$. The {\em quotient graph} $\G_N$ of $\G$ relative to
$N$ is defined as the graph with vertices the orbits of $N$ on $V(\G)$ and with two orbits
adjacent if there is an edge in $\G$ between those two orbits. The theory of quotient graph
is widely used to investigate symmetric graphs.
Let $\G$ be a symmetric graph and let $N$ be a normal subgroup of $\Aut(\G)$. If $\G$ and $\G_N$ have the
same valency, then the graph $\G$ is said to be a {\em normal cover} of $\G_N$ and the graph $\G_N$ is called
a {\em normal quotient} of $\G$. In this case, $N$ is semiregular on $V(\G)$.
In particular, a graph $\G$ is called {\em basic} if it has no proper normal
quotient.

%It is well known that a locally primitive graph is basic if and only if it has no nontrivial normal subgroup having more than two orbits.

To study a symmetric graph $\G$, there is a natural strategy,
which contains the following two steps:
\begin{itemize}
\item Step~1: Investigate normal quotient graph
$\G_N$ for some normal subgroup $N$ of $\Aut(\G)$;
\item Step~2: Reconstruct
the original graph $\G$ from the normal quotient $\G_N$ by using covering
techniques.
\end{itemize}
For Step~1, it is usually done by taking the normal subgroup $N$
as large as possible and then the graph $\G$ is reduced a basic graph.
In this paper, we only consider Step~1, that is, basic graphs;
for Step~2, one may refer~\cite{AHK,FK,IKM}.

Symmetric graphs of order a small number  times a prime power
have been received a lot of attention. Such graphs not only provide
plenty of significant examples of graphs (see~\cite{IK,PLY,ZF}) but also have been succeeded in investigating many graphs with different symmetric properties (see~\cite{FLZ}).
Let $p$ be a prime.
Recently, Morgan et al.~\cite{MSV} proved that
there are only finitely many $d$-valent 2-arc-transitive graphs of order $kp^n$
for given integers $n$ and $k$ if $d$ is large enough. Moreover, they pointed that for classifying 2-arc-transitive graphs along these lines, the most interesting
case is when the valency is small.

In the literature, there are lots of works in classifying symmetric graphs of order $kp^n$ and small valency $d$, especially with $n$
and $k$ small and $d=3$, 4, or 5; see~\cite{FK6,FZL, GZF,GHS, PLY, PLHL, ZF} for example.
In particular, cubic or pentavalent symmetric basic graphs of order $2p^n$
for variable prime $p$ and integer $n$ were classified completely in~\cite{FK6,FZL},
and tetravalent 2-arc-transitive basic graphs of order $2p^n$
were classified in~\cite{ZF}. As applications of above works,
symmetric graphs of order $2p^2$ and valency 3, 4 or 5 were determined
(for valency 5, one may also refer~\cite{PLY}).
Cubic and tetravalent symmetric graphs of order $kp^n$
for some small integer $k$ with $k>2$ were also well studied (see~\cite{CLP, LWX,PLHL}).
This paper concerns pentavalent symmetric graphs of order $kp^n$ with $k>2$, and the following theorem is the main result.

\begin{theorem}\label{theo=2qp^n}
Let $p>q$ be two primes and $n$ a non-negative integer. Let $\G$ be a connected pentavalent symmetric graph of order $2qp^n$. Then $\G$ is a normal cover of  one graph listed in Table~\ref{table=2}. Moreover, all graphs in Table~\ref{table=2} are basic except for $K_{6,6}-6K_2$ and $\mathbf{I}_{12}$.

\begin{table}[h]

\begin{center}
\begin{tabular}{|c|c|c|c|c|c|}

\hline
$\Sigma$ & $(p,q)$ & $\Aut(\Sigma)$ & $\Sigma$ & $(p,q)$ & $\Aut(\Sigma)$ \\
\hline
$K_6$ & $q=3$, $p>3$ & $S_6$& $K_{5,5}$ & $q=5$, $p>5$ & $(S_5\times S_5)\rtimes\mz_2$\\
\hline
$\mathcal{CD}_q$ & $5\di (q-1)$, $p>q$& $D_q\rtimes\mz_5$& $K_{6,6}-6K_2$ & $(3,2)$
&$S_6\times\mz_2$ \\
\hline
 $\mathbf{I}_{12}$ & $(3,2)$  & $A_5\times\mz_2$ & $\mathcal{G}_{36}$ & $(3,2)$
 &$\Aut(A_6)$ \\
\hline
$\mathcal{G}_{42}$ & $(7,3)$ & $\Aut(\PSL(3,4))$& $\mathcal{G}_{66}$ & $(11,3)$ & $\PGL(2,11)$\\
\hline
 $\mathcal{G}_{114}$ & $(19,3)$ & $\PGL(2,19)$& $\mathcal{G}_{170}$ & $(17,5)$ &$\Aut(\PSp(4,4))$\\
\hline
 $\mathcal{G}_{406}$ & $(29,7)$ & $\PGL(2,29)$ &$\mathcal{G}_{574}$ & $(41,7)$ & $\PSL(2,41)$\\
\hline
$\mathcal{G}_{3422}$ & $(59,29)$ & $\PGL(2,59)$ & $\mathcal{G}_{3782}$ & $(61,31)$ & $\PGL(2,61)$  \\
\hline
\end{tabular}
\end{center}
\vskip -0.5cm
\caption{{\small Normal quotients of symmetric pentavalent graphs of order $2pq^n$}}\label{table=2}
\end{table}

\end{theorem}

All graphs in Table~\ref{table=2} will be introduced in Section~2.
For pentavalent basic graphs of order $2qp^n$, the graphs graphs with $q=p$ have been determined in~\cite{FZL}.
It is known from Conder et al.~\cite{CLP} that there are only
finitely many connected pentavalent 2-arc-transitive graphs of order
$kp$ or $kp^2$ for a given integer $k$.
Since the graph $\mathcal{CD}_q$ is not 2-arc-transitive,
by Theorem~\ref{theo=2qp^n}, one may easily see that
for any given integer $n$, there are finitely many pentavalent 2-arc-transitive graph
of order $kp^n$ if $k=2q$ with prime $q$ such that $5<q<p$.
This does not hold for $k=2p$, see~\cite[Theorem~5.3]{FZL}.
Moreover, as an application of Theorem~\ref{theo=2qp^n},
connected pentavalent symmetric graphs of order $4p^n$ with $2\leq n\leq 4$ or $6p^2$ are classified completely.

\section{Preliminaries}

In this section, we describe some preliminary results which will be used later.
Denote by $\mz_n$, $D_n$, $A_n$ and $S_n$ the cyclic group of order $n$, the dihedral group of order $2n$, the alternating group and
the symmetric group of degree $n$, respectively.

%Let $k$ be a given integer. Conder et al.~\cite{CLP} proved that there are only
%finitely many connected 3-valent 2-arc-transitive graphs of order
%$kp$ for some prime $p$, and there are finitely many connected $d$-valent
%2-arc-transitive graphs of order $kp$ or $kp^2$ if $d\geq 4$.

\subsection{Group theory}

Let $G$ and $E$ be two groups. We call an extension $E$ of $G$ by $N$ a
{\em central extension} of $G$ if $E$ has a central subgroup $N$ such that $E/N\cong G$,
and if further $E$ is perfect, that is, the derived group $E'=E$, we
call $E$ a {\em covering group} of $G$. Schur proved that for every non-abelian simple group $G$
there is a unique maximal covering group $M$ such that every covering group of $G$ is a factor
group of $M$ (see \cite[V \S23]{Huppert}). This group $M$ is called the {\em full covering group} of $G$, and the center of $M$
is the {\em Schur multiplier} of $G$, denoted by $\Mult(G)$.
The following proposition can be obtained from~\cite[Lemma~2.11]{PLHL}.

\begin{prop}\label{lem=m/n}
Let $M$ be a finite group and $N$ a normal subgroup of $M$ of order $p$ or $p^2$, where $p$ is a prime.
If $M/N$ is a non-abelian simple group, then $M=M'N$ and either $M$
is a covering group of $M/N$ and $N\lesssim \Mult(M/N)$, or $p\di |M:M'|$.
\end{prop}

Let $G$ be a group, and let $\pi(G)$ be the set of distinct prime divisors
of $|G|$. The group $G$ is called a {\em $K_n$-group}
if $|\pi(G)|=n$. The following lemma is based on the classification
of finite non-abelian simple $K_n$-groups with $n\leq 5$.

%By~\cite[Theorem~I]{HL}, there exist 8 non-abelian 3-prime factor
%sporadic simple groups and 31 non-abelian 4-prime factor sporadic
%simple groups except for $\PSL(2,p^n)$, where $p$ is a prime and $n$
%is a positive integer. For a 4-prime factor simple group
%$\PSL(2,p^n)$ having a prime divisor 5, by~\cite[Theorem~3.2,
%Lemmas~3.5~(2) and 3.4~(2)]{HL}, we have either $p^n\in
%\{2^4,3^4,5^2,7^2\}$ or $p\geq11$ and $n=1$.
%By~\cite[Proposition~2.3]{GZF}, one can obtain the following
%proposition (also see \cite[pp.134-136]{Gorenstein}).

\begin{lem}\label{prop=simple}
Let $p>q$ be two primes, and let $G$ be a non-abelian simple group of order $2^i\cdot 3^j\cdot 5\cdot q^k\cdot p^{\ell}$ with integers
$1\leq i\leq10$, $0\leq j\leq 2$, $0\leq k\leq 1$,
and $\ell\geq 0$. Then one of the following holds.

\begin{enumerate}

\itemsep -1pt
\item [\rm (1)] $|\pi(G)|=3$ and $G\cong A_5$, $A_6$ or $\PSU(4,2)$
with order $2^2\cdot 3\cdot 5$, $2^3\cdot 3^2\cdot 5$, or $2^6\cdot 3^4\cdot 5$,
respectively.
\item [\rm (2)] $|\pi(G)|=4$ and $G$ is isomorphic to one group listed in Table~{\rm\ref{table=1}}.
\item [\rm (3)] $|\pi(G)|=5$ and $k=\ell=1$.

\end{enumerate}

\begin{table}[h]

\begin{center}
\begin{tabular}{|l|l|l|l|l|l|}

\hline

$G$ & Order & $G$ & Order & $G$ & Order\\
\hline
$A_7$ & $2^3\cdot 3^2\cdot 5\cdot 7$ & ${\rm M}_{12}$ & $2^6\cdot 3^3\cdot 5\cdot 11$ &$\PSL(2,5^2)$ &$2^3\cdot3\cdot5^2\cdot13$ \\
$A_8$ & $2^6\cdot 3^2\cdot 5\cdot 7$ &$\PSU(3,4)$& $2^6\cdot3\cdot5^2\cdot13$&
$\PSL(2,7^2)$ & $2^4\cdot 3\cdot 5^2\cdot 7^2$\\
$\PSL(3,4)$ &$2^6\cdot3^2\cdot5\cdot7$ &$Sp(4,4)$& $2^8\cdot3^2\cdot5^2\cdot13$ &$\PSL(2,11)$ &$2^2\cdot3\cdot5\cdot11$\\
$\PSp(4,7)$& $2^8\cdot 3^2\cdot 5^2\cdot7^4$ & $Sz(8)$ & $2^6\cdot 5\cdot7\cdot 13$&
$\PSL(2,19)$ & $2^2\cdot 3^2\cdot 5\cdot 19$\\
${\rm M}_{11}$ & $2^4\cdot 3^2\cdot 5\cdot 11$ & $\PSL(2,2^4)$& $2^4\cdot 3\cdot 5\cdot 17$&
$\PSL(2,31)$ & $2^5\cdot 3\cdot 5\cdot 31$\\
\hline
\end{tabular}
\end{center}
\vskip -0.5cm
\caption{{\small Some simple $K_4$-groups}}\label{table=1}
\end{table}

\end{lem}

\f {\bf Proof:} Clearly, $|\pi(G)|=3$, 4 or 5.
For $|\pi(G)|=3$, by~\cite[Theorem~I]{HL}, $G$ is isomorphic to one group
listed in~\cite[Table 1]{HL}, and since $5\in \pi(G)$,
we have $G\cong A_5$, $A_6$, or $\PSU(4,2)$.

Assume that $|\pi(G)|=4$. Then $p\geq7$.
Again by~\cite[Theorem~I]{HL}, $G$ is isomorphic to one group
listed in~\cite[Table 2]{HL} or a non-abelian simple $K_4$-group $\PSL(2,r^n)$
for some prime $r$ and positive integer $n$.
Note that $|G|=2^i\cdot 3^j\cdot 5\cdot q^k\cdot p^{\ell}$
with  $0\leq j\leq 2$ and $0\leq k\leq 1$. Since $p\geq7$, we have $3^4\nmid |G|$
and $5^3\nmid|G|$. In particular, $3^3\cdot 5^2\nmid |G|$, and if $p,q\geq 7$, then $5^2\nmid |G|$. For the sporadic simple groups, by~\cite[Table 2]{HL} we have $\G\cong A_7$, $A_8$,
$\PSL(3,4)$, $\PSp(4,7)$,
${\rm M}_{11}$, ${\rm M}_{12}$, $\PSU(3,4)$, $Sp(4,4)$, or $Sz(8)$.
For the simple $K_4$-group $\PSL(2,r^n)$,
since $5\di |G|$ and $3^4\nmid|G|$, by~\cite[Theorem~3.2,
Lemmas~3.5~(2) and 3.4~(2)]{HL} either $r^n\in
\{2^4,5^2,7^2\}$ or $r\geq11$ and $n=1$.

Let $G\cong \PSL(2,r)$ with $r\geq11$. Then
$|G|=\frac{1}{2}r(r-1)(r+1)$. Since either $r-1$ or $r+1$ has a prime divisor 3,
the group $G$ is a $\{2,3,5,r\}$-group,
and since $|G|=2^{i}\cdot 3^j\cdot 5\cdot q^k\cdot p^{\ell}$ with $1\leq i\leq10$, $0\leq j\leq 2$, $0\leq k\leq 1$ and $p\geq7$, we
have $\frac{1}{2}(r-1)(r+1)\di 2^{11}\cdot 3^3\cdot 5^2$.
Note that $(\frac{r-1}{2},\frac{r+1}{2})=1$ and $r\geq11$.
If $2\di \frac{r-1}{2}$, then $\frac{r+1}{2}\di 3^3\cdot 5^2$,
and thus $r\in\{17,29,53,89,149,269,449\}$.
However, for each $r\in\{17,29,53,89,149,269,449\}$,
the group $\PSL(2,r)$ is not a simple $K_4$-group
such that $5\di |\PSL(2,r)|$ (this can be also checked by MAGMA~\cite{BCP} easily).
If $2\di \frac{r+1}{2}$, then $\frac{r-1}{2}\di 3^3\cdot 5^2$,
and similarly, we have $r\in\{11,19,31\}$.
Hence $G\cong\PSL(2,11)$, $\PSL(2,19)$ or $\PSL(2,31)$.

Assume that $|\pi(G)|=5$. Then $p>q\geq 7$
and $|G|=2^i\cdot 3^j\cdot 5\cdot q\cdot p^{\ell}$ with integers
$1\leq i\leq10$, $1\leq j\leq 2$, and $\ell\geq 1$, implying that
\begin{equation}\label{eq=1}
2^{11}\nmid |G|,~3^3\nmid |G|,~5^2\nmid|G|,~{\rm and}~q^2\nmid |G|.
\end{equation}

To finishe the proof, it suffices to prove that $\ell=1$.
By~\cite[Theorem~A]{JI}, $G$ is isomorphic to one of 30 sporadic simple groups,
$\PSL(2,r)$, $\PSL(3,r)$, $\PSU(3,r)$, $O_5(r)$, $Sz(2^{2m+1})$ with $m\geq3$,
or $R(3^{2m+1})$ with $m\geq1$, where $r$ is a prime power.

 For the 30 sporadic simple groups,
by Eq~(\ref{eq=1}) we have $G\cong {\rm M}_{22}$ or $\PSL(5,2)$
with order $2^7\cdot 3^2\cdot 5\cdot 7\cdot 11$ or
$2^{10}\cdot 3^2\cdot 5\cdot 7\cdot 31$, respectively, implying that
$p=11$ or 31 and $\ell=1$.

Suppose that $G\cong Sz(2^{2m+1})$ with $m\geq3$ or $R(3^{2m+1})$ with $m\geq1$.
Then $|G|=2^{4m+2}(2^{4m+2}+1)(2^{2m+1}-1)$
or $3^{6m+3}(3^{6m+3}+1)(3^{2m+1}-1)$ (see \cite[pp.134-136]{Gorenstein}), respectively,
implying that $2^{11}\di |G|$
or $3^3\di |G|$, contrary to Eq~(\ref{eq=1}).
Hence $G\ncong Sz(2^{2m+1})$ or $R(3^{2m+1})$.

Suppose that $G\cong O_5(r)$. Then $|G|=\frac{1}{2}r^4(r^4-1)(r^3-1)(r^2-1)$.
By Eq~(\ref{eq=1}), $r$ is a 2-power or $p$-power.
If $r$ is a 2-power, then $r=2$ or $2^2$
because $2^{11}\nmid |G|$ (see Eq~(\ref{eq=1})).
Since $G$ is a $K_5$-group, we have $r=2^2$,
and so $3^4\di |G|$, contrary to Eq~(\ref{eq=1}).
If $r$ is a $p$-power, then $\frac{1}{2}(r^{4}-1)(r^{3}-1)(r^{2}-1)\di 2^{10}\cdot 3^2\cdot 5\cdot q$, that is, $(r^{2}+1)(r^{3}-1)(r^{2}-1)^2\di 2^{11}\cdot 3^2\cdot 5\cdot q$.
Since $q^2\nmid |G|$ and $5^2\nmid |G|$ by Eq~(\ref{eq=1}),
we have $q\nmid (r^2-1)$ and $5\nmid (r^2-1)$,
forcing that $(r^{2}-1)^2\di 2^{11}\cdot 3^2$
and $(r^{2}-1)\di 2^{5}\cdot 3$.
This is impossible because $r$ is a $p$-power with $p>7$.
Hence $G\ncong O_5(r)$.

Suppose that $G\cong \PSU(3,r)$. Then $|G|=\frac{1}{(3,r+1)}r^3(r-1)(r+1)^2(r^2-r+1)$,
and by Eq~(\ref{eq=1}), $r$ is a 2-power or a $p$-power.
If $r$ is a 2-power, then $r=2$, $2^2$ or $2^3$,
and the corresponding groups $\PSU(3,r)$
is not a simple $K_5$-group, a contradiction.
If $r$ is a $p$-power, then $\frac{1}{(3,r+1)}(r-1)(r+1)^2(r^2-r+1)\di
2^{10}\cdot 3^2\cdot 5\cdot q$,
implying that $(r+1)^2\di 2^{10}\cdot 3^3$
and $(r+1)\di 2^5\cdot 3$.
Since $r$ is a $p$-power with $p>7$,
we have $r=11$, 23, 31 or 47,
and since $G\cong \PSU(3,r)$ is a $K_5$-group
and $5\di |G|$,
we have $G\cong\PSU(3,11)$ and $|G|=2^5\cdot 3^2\cdot 5\cdot11^3\cdot37$,
yielding that $q=37$, contrary to that $q<p$.
For $\PSL(3,r)$,
we have $|\PSL(3,r)|=\frac{1}{(3,r-1)}r^3(r-1)^2(r+1)(r^2+r+1)$,
and similarly, one may check that $G\ncong\PSL(3,r)$.

Finally, suppose that $G\cong \PSL(2,r)$. Then $|G|=\frac{1}{(2,r-1)}r(r-1)(r+1)$.
Since $p>q\geq7$ and $p\di |G|$, we have $r\neq 3$, $3^2$, 5 or $q$.
By Eq~(\ref{eq=1}), $r=2^s$ with $s\leq 10$
or $r=p^{\ell}$.
For $r=2^s$, by checking the order of $\PSL(2,r)$,
we have $G\cong\PSL(2,2^6)$ or $\PSL(2,2^8)$
of order $2^6\cdot 3^2\cdot5\cdot7\cdot13$
or $2^8\cdot 3\cdot5\cdot17\cdot257$,
implying that $p=13$ or 257 and $\ell=1$, as required.
For $r=p^{\ell}$, we have
$\frac{1}{2}(r-1)(r+1)\di  2^{10}\cdot 3^2\cdot 5\cdot q$.
Since $(\frac{r-1}{2},\frac{r+1}{2})=1$,
either $(r-1)\di 2^{10}\cdot 3^2\cdot 5$
or $(r+1)\di 2^{10}\cdot 3^2\cdot 5$.
Recall that $p>7$.
If $\ell\geq 2$, then
$r=11^2$, $17^2$, $19^2$, or $31^2$
(this can be checked easily by MAGMA~\cite{BCP}),
and since $G\cong\PSL(2,r)$ is a $K_5$-group,
$r=11^2$, $17^2$ or $19^2$
and $|G|=2^3\cdot 3\cdot5\cdot 11^2\cdot 61$, $2^5\cdot 3^2\cdot5\cdot 17^2\cdot 29$
or $2^3\cdot 3^2\cdot5\cdot 19^2\cdot 181$, contrary to that
$p$ is the largest prime divisor of $G$.
If follows that $\ell=1$.
\hfill\qed

%and since $(\frac{r-1}{2},\frac{r+1}{2})=1$,
%$r-1$ has at most three distinct prime divisors.
%Recall that $p>7$. If $\ell\geq 2$ is even, then~\cite[Lemma~1.4]{HL}
%implies that $r=p^{\ell}=5^4$ or $7^4$, a contradiction.
%If $\ell\geq 3$ is odd, then
%$(r-1)(r+1)=(p-1)(p^{\ell-1}+p^{\ell-2}+\cdots +1)(p+1)(p^{\ell-1}-p^{\ell-2}+\cdots +p^2-p+1)$.
%In particular, both $(p^{\ell-1}+p^{\ell-2}+\cdots +1)$ and $(p^{\ell-1}-p^{\ell-2}+\cdots + p^2-p+1)$ are odd and they are relatively prime.
%It forces that $(p^{\ell-1}+p^{\ell-2}+\cdots +1)\di 3^2\cdot 5$
%or $(p^{\ell-1}-p^{\ell-2}+\cdots + p^2-p+1)\di 3^2\cdot 5$,
%which is impossible because $p>7$. Hence $\ell=1$.
%

\subsection{Pentavalent symmetric graphs}

The following proposition is about the stabilizers of connected pentavalent symmetric graphs,
coming from~\cite[Theorem 1.1]{GF}.

\begin{prop}\label{prop=stabilizer}
Let $\G$ be a connected pentavalent $(G,s)$-transitive graph for $G\leq \Aut(\G)$ and $s\geq1$.
Let $v\in V(\G)$. Then $s\leq 5$, $|G_v|\di 2^9\cdot 3^2\cdot5$ and one of the following holds.
\begin{enumerate}

\itemsep -1pt
\item [\rm (1)] For $s=1$, $G_v\cong \mathbb{Z}_5$, $D_{5}$ or $D_{10}$. In particular, $|G_v|=5$, $2\cdot 5$ or $2^2\cdot 5$.
\item [\rm (2)] For $s=2$, $G_v\cong F_{20}$, $F_{20}\times \mathbb{Z}_2$, $A_5$ or $S_5$,  where $F_{20}$ is the Frobenius group of order $20$. In particular, $|G_v|=2^2\cdot 5$, $2^3\cdot 5$, $2^2\cdot 3 \cdot 5$ or $2^3\cdot 3 \cdot 5$.
\item [\rm (3)] For $s=3$, $G_v\cong F_{20}\times \mathbb{Z}_4$, $A_4\times A_5$, $S_4\times S_5$ or $(A_4\times A_5)\rtimes \mathbb{Z}_2$ with $A_4\rtimes \mathbb{Z}_2=S_4$ and $A_5\rtimes \mathbb{Z}_2=S_5$. In particular, $|G_v|=2^4\cdot 5$, $2^4\cdot 3^2\cdot 5$, $2^6\cdot 3^2\cdot 5$ or $2^5\cdot 3^2\cdot 5$.
\item [\rm (4)] For $s=4$, $G_v\cong {\rm ASL}(2,4)$, ${\rm AGL}(2,4)$, ${\rm A\Sigma L}(2,4)$ or ${\rm A\Gamma L}(2,4)$. In particular, $|G_v|=2^6\cdot3 \cdot 5$, $2^6\cdot 3^2\cdot 5$, $2^7\cdot 3\cdot 5$ or $2^7\cdot 3^2\cdot 5$.
\item [\rm (5)] For $s=5$, $G_v\cong \mathbb{Z}^6_2\rtimes {\rm \Gamma L}(2,4)$ and $|G_v|=2^9\cdot 3^2\cdot 5$.
\end{enumerate}
\end{prop}

By~\cite[Lemma~5.1]{FZL}, we have the following proposition.

\begin{prop}\label{prop=minimal normal subgroup}
Let $p$ be a prime and $n\geq2$ an integer. Let $\Gamma$ be a connected pentavalent $G$-arc-transitive graph of
order $2p^n$ with $G\leq \Aut(\G)$. Then every minimal
normal subgroup of $G$ is an elementary abelian $p$-group.
\end{prop}

In view of \cite[Theorem 9]{Lorimer}, we have the following.

\begin{prop} \label{prop=atlesst3orbits}
Let $\Gamma$ be a connected pentavalent $G$-arc-transitive graph, and let $N\unlhd G$.
If $N$ has at least three orbits, then it is semiregular on $V(\G)$ and the quotient graph $\G_N$ is a connected pentavalent $G/N$-arc-transitive graph.
\end{prop}

%Next, we introduce some pentavalent symmetric graphs
%constructed as coset graphs or Cayley graphs,
%which will be used later.
Let $G$ be a finite
group and $H\leq G$. Denote by $D$ a union of some double cosets of
$H$ in $G$ such that $D^{-1}=D$. A {\em coset graph}
$\G=\Cos(G,H,D)$ on $G$ with respect to $H$ and $D$ is defined to
have vertex set $V(\G)=[G:H]$, the set of right cosets of $H$ in
$G$, and the edge set $E(\G)=\{\{Hg,Hdg\}~|~g\in G, d\in D\}$.
The following two coset
graphs were introduced in~\cite[Examples~3.5 and 3.6]{HFL}.

\begin{exam} {\rm \label{exam=G42}
Let $T=\PSL(3, 4)$. Then $T$ has two conjugate classes of maximal
subgroups isomorphic to $\mz_2^4\rtimes A_5$, which are fused by an automorphism $g$ of $T$ of order 2. Let $H=\mz_2^4\rtimes A_5$ be a subgroup of $T$ and $G=\lg T, g\rg$.
By~\cite[Example~3.5]{HFL}, the coset graph $\Cos(G,H,HgH)$, denoted by $\mathcal{G}_{42}$,
is a connected pentavalent 4-transitive graph of order $42$ with $\Aut(\mathcal{G}_{42})\cong\Aut(T)$.
}
\end{exam}

\begin{exam} {\rm \label{exam=G170}
Let $T=\PSp(4, 4)$. Then $\Inn(T)\cong T$ has two conjugate classes of maximal
subgroups isomorphic to $\mz_2^6\rtimes (\mz_3\times A_5)$, which are fused
in $\Aut(T)$. Let $f$ be a field automorphism, and let
$g\in\Aut(T)/(\Inn(T).\lg f\rg)$.
Then $f$ has order 2, $\Aut(T)/(\Inn(T).\lg f\rg)\cong\mz_2$
and $g^2\in\lg f\rg$.
Let $S\cong \mz_2^6\rtimes (\mz_3\times A_5)$
be a subgroup of $\Inn(T)$ such that $f$ normalizes $S$ but $g$ does not.
Set $G=\lg \Inn(T), f, g\rg$, and $H=\lg S, f\rg$.
By~\cite[Example~3.6]{HFL}, the coset graph $\Cos(G, H, HgH)$,
denoted by $\mathcal{G}_{170}$, is a connected pentavalent 5-transitive graph
and $\Aut(\mathcal{G}_{170})\cong\Aut(T)$.

}
\end{exam}

The following graph $\mathcal{G}_{36}$ comes from~\cite[Section~3]{GZF},
and the graph $\mathcal{G}_{108}$ is a normal cover of $\mathcal{G}_{36}$,
coming from~\cite{P}.

\begin{exam}{\rm  \label{exam=G36}
Let $G\cong A_6$ and $P$ a Sylow 5-subgroup of $G$. Let $H=N_G(P)$.
From Atlas~\cite[pp.~4]{Atlas}, we have $H\cong D_5$ and for any involution $x\in H$,
the centralizer of $x$ in $G$ is isomorphic to $D_4$, that is, $C_G(x)\cong D_4$. Let
$g$ be an element of order $4$ in $C_G(x)$.
Denoted by $\mathcal{G}_{36}$ the coset graph $\Cos(G,H,HgH)$.
By~\cite[Theorem~4.1]{GZF}, $\mathcal{G}_{36}$ is a connected pentavalent $2$-transitive graph of order $36$ and $\Aut(\mathcal{G}_{36})\cong\Aut(A_6)$. Moreover, $\mathcal{G}_{36}$
is the unique connected pentavalent symmetric graph of order $36$
up to isomorphic.

}
\end{exam}

%Next, we introduce a new pentavalent symmetric graph, which is a normal cover
%of $\mathcal{G}_{36}$ and has order 108.

\begin{exam} {\rm \label{exam=G108}
By Atlas~\cite[pp.~4]{Atlas}, we have $\Mult(A_6)\cong\mz_6$. Let $M$ be the full covering group
of $A_6$. Then $M$ has a central subgroup $\mz_2$. Set $T=M/\mz_2$.
Then $T$ is a covering group of $A_6$ with center $\mz_3$.
Let $H$ be a subgroup of $T$ of order $10$ and $x$ an involution in $H$. By Example~\ref{exam=G36}, we have $C_T(x)\cong D_4\times \mz_3$.
Take an element $g$ of order $4$ in $C_T(x)$. Then $g^2=x\in H$ and $\lg H,g\rg=T$.
Denote by $\mathcal{G}_{108}$ the coset graph $\Cos(T,H,HgH)$.
By {\rm MAGMA~\cite{BCP}}, $\mathcal{G}_{108}$ is a connected pentavalent $2$-transitive graph of order $108$ and $\Aut(\mathcal{G}_{108})\cong \mz_3.\Aut(A_6)$.}
\end{exam}

The following six graphs were constructed in~\cite[Section~3]{HFL}.

\begin{exam}{\rm \label{exam=G66}
Let $G$ be a group and $H$ a subgroup of $G$ of index $n$, as listed in Table~\ref{table=3}. Then $G$ has an involution $g$ such that $|HgH|/|H|=5$ and $\lg H, g\rg=G$, and the coset graph $\Cos(G, H, HgH)$,
denoted by $\mathcal{G}_{n}$, is a connected pentavalent $s$-transitive graph
of order $n$.

}
\end{exam}

\begin{table}[h]

\begin{center}
\begin{tabular}{|l|l|l|l|l|}

\hline
$G$& $H$ & $\mathcal{G}_{n}$ & $s$ &$\Aut(\mathcal{G}_{n})$\\
\hline
$\PSL(2,11)$ & $D_5$ & $\mathcal{G}_{66}$ & $1$ &$\PGL(2,11)$ \\
 \hline
 $\PGL(2,19)$ & $A_5$ & $\mathcal{G}_{114}$ & $2$ & $\PGL(2,19)$\\
\hline
$\PGL(2,29)$ & $A_5$ & $\mathcal{G}_{406}$ & $2$ & $\PGL(2,29)$ \\
\hline
$\PGL(2,59)$& $A_5$ & $\mathcal{G}_{3422}$ & $2$ & $\PGL(2,59)$\\
\hline
$\PGL(2,61)$& $A_5$ & $\mathcal{G}_{3782}$ & $2$ & $\PGL(2,61)$ \\
\hline
$\PSL(2,41)$& $A_5$ & $\mathcal{G}_{574}$ & $2$ & $\PSL(2,41)$\\
\hline
\end{tabular}
\end{center}
\vskip -0.5cm
\caption{{\small Some connected pentavalent symmetric graphs}}\label{table=3}
\end{table}

Let $G$ be a group, and let $S$ be a generated subset of $G$ with $1\notin S$ and $S^{-1}=S$.
Clearly, the coset graph $\G=\Cos(G,1,S)$ is a connected undirected simple graph,
which is called a {\em Cayley graph} and denoted by $\Cay(G,S)$.
The following infinite
family of Cayley graphs was first constructed in~\cite{KKO}.

%Let $R(G)$ be the right regular representation of $G$,
%the acting group of $G$ by right multiplication. Then $R(G)\leq \Aut(\G)$.
%By Godsil~\cite{G} or Xu~\cite{X},
%$N_{{\rm Aut}(\G)}(R(G))=R(G)\rtimes\Aut(G,S)$, where $\Aut(G,S)=\{\a\in \Aut(G) \ |\ S^\a=S\}$. A Cayley graph
%$\G=\Cay(G,S)$ is said to be {\em normal} if $R(G)$ is normal in $\Aut(\G)$, and in
%this case, $\Aut(\G)=R(G)\rtimes\Aut(G,S)$.

%For an abelian group $H$,
%the {\em generalized dihedral group} ${\rm Dih}(H)$ is the semidirect product $H\rtimes \mz_2$,
%where the unique involution in $\mz_2$ maps each element of $H$ to its inverse. In
%particular, if $H$ is cyclic, ${\rm Dih}(H)$ is a dihedral
%group.

\begin{exam}{\rm \label{exam=n}
Let $m>1$ be an integer such that $x^4+x^3+x^2+x+1=0$ has a
solution $r$ in $\mz_m$. Then $m=5$, 11 or $m\geq31$. Let
$$\mathcal{CD}_{m}=\Cay(D_{m},\{b,ab,a^{r+1}b,a^{r^2+r+1}b,a^{r^3+r^2+r+1}b\})$$
be a Cayley graph on the dihedral group
$D_{m}=\lg a,b~|~a^{n}=b^2=1, a^b=a^{-1}\rg$.
For $m=5$ or 11, by~\cite{Cheng},
$\Aut(\mathcal{CD}_m)\cong (S_5\times S_5)\rtimes\mz_2$
or $\PGL(2,11)$, respectively. In particular, $\mathcal{CD}_{5}\cong K_{5,5}$,
the complete bipartite graph of order 10.
For $m\geq31$, by~\cite[Theorem~B and Proposition~4.1]{KKO},
$\Aut(\mathcal{CD}_{m})\cong D_{m}\rtimes\mz_5$.

}\end{exam}

Denote by $\mathbf{I}_{12}$ the Icosahedron graph
and by $K_{6,6}-6K_2$ the complete bipartite graph
of order 12 minus a one-factor. By~\cite[Proposition 3.2]{GZF},
$\Aut(\mathbf{I}_{12})\cong A_5\times \mz_2$
and $\Aut(K_{6,6}-6K_2)\cong S_6\times \mz_2$.
Based on~\cite{Cheng} and~\cite{HFL}, we have the following proposition.

\begin{prop}\label{prop=2p}
Let $\G$ be a connected  pentavalent symmetric graph. Then the following hold.
\begin{enumerate}
  \item[{\rm (1)}] If $|V(\G)|=2p$ for a prime $p$, then $\G\cong K_6$,
   $K_{5,5}$ or $\mathcal{CD}_{p}$ with $5\ |\ (p-1)$.
  \item [\rm (2)] If $|V(\G)|=2pq$ for two primes $p>q$, then
$\G\cong K_{6,6}-6K_2$, $\mathbf{I}_{12}$, $\mathcal{G}_{42}$,
$\mathcal{G}_{170}$, $\mathcal{CD}_{pq}$ with $5\di (p-1)$
and $q=5$ or $5\di (q-1)$, or one graph listed in Table~\ref{table=3}.
\end{enumerate}

\end{prop}

\section{Proof of Theorem~\ref{theo=2qp^n}}

In this section, we aim to prove Theorem~\ref{theo=2qp^n}.
First, we need the following lemma.

\begin{lem}\label{lem=4pn-minimal normal subgroup}
Let $p>q$ be two primes, and let $\G$ be a connected pentavalent $G$-arc-transitive graph of order $2qp^n$ with $n\geq2$.
If $\G\ncong \mathcal{G}_{36}$, then $G$ has a minimal normal subgroup that is an elementary abelian $p$-group.
\end{lem}

\f {\bf Proof:} Let $N$ be a minimal normal subgroup of $G$. Then $N=T^s$ for a positive integer $s$ and a simple group $T$.

Suppose that $N$ is insolvable. Then $T$ is non-abelian. Let $v$ be a vertex of $\G$.
Since $\G$ has the prime valency $5$, $G_{v}$ is primitive on the neighborhood of $v$ in $\G$,
and since $N_v\unlhd G_v$, either $N_v=1$ or $5\di |N_v|$.
If $N_v=1$, then $|N|\di 2qp^n$. Note that $p>q$ implies that $2qp^n=2^2p^n$, or $2qp^n$ is twice an odd integer. It follows that $N$ is solvable, a contradiction. Hence $5\di |N_v|$.
By Proposition~\ref{prop=atlesst3orbits}, $N$ has at most two orbits on $V(\G)$,
forcing that $qp^n\di |v^N|$ and $5qp^n\di |N|$, where $v^{N}$ is the orbit of $N$ on $V(\G)$ containing the vertex $v$.

Since $|G_v|\di 2^9\cdot3^2\cdot5$ by Proposition~\ref{prop=stabilizer},
we have $|G|\di 2^{10}\cdot 3^2\cdot 5\cdot q\cdot p^n$,
and since $N\leq G$ and $5qp^n\di |N|$, we have
\begin{equation}\label{eq=2}
|N|=2^i\cdot 3^j\cdot 5\cdot q\cdot p^n~{\rm with}~0\leq i\leq 10,~0\leq j\leq 2,~{\rm and}~n\geq2.
\end{equation}
In particular, $|\pi(N)|=3$, 4 or 5, where $\pi(N)$ is the set of distinct prime factors of $|N|$. Recalling that $N=T^s$, we have $5\di |T|$, $q\di |T|$ and $p\di |T|$
(note that $p$ or $q$ may be equal to 5).
We consider the following three cases depending on $|\pi(N)|$.

\medskip
\f {\bf Case~1:} $|\pi(N)|=3$.

In this case, $T$ is a non-abelian simple
$K_3$-group, and since $5\di |T|$,
we have $T\cong A_5$, $A_6$ or $\PSU(4,2)$ by Lemma~\ref{prop=simple}~(1).
It implies that $p,q\in\{2,3,5\}$, and since $p>q$, we have $p=5$ or 3.

If $p=5$, then $q=2$ or 3. By Eq~(\ref{eq=2}),
we have $|N|=2^i\cdot 3^j\cdot 5^{n+1}\cdot q$ with $0\leq j\leq 2$,
implying that $N=T^{n+1}$ with $n+1\geq3$ and $3^4\nmid |N|$.
Hence $N\cong A_5\times A_5\times A_5$ with $q=3$ and $n=2$,
yielding that $|V(\G)|=6\cdot 5^2$.
Clearly, $A_5\cong T\lhd N$ has more
than two orbits on $V(\G)$ and $T$ is not semiregular because $|T|=60$.
Recall that $N$ has at most two orbits on $V(\G)$ and $5\di |N_v|$.
If $N$ is transitive on $V(\G)$, then $N$ is arc-transitive
and Proposition~\ref{prop=atlesst3orbits} implies that $T$ is semiregular, a contradiction.
If $N$ has exactly two orbits on $V(\G)$, then $\G$ is bipartite with the
orbits of $N$ as its bipartite sets. Since $|T|=60$, $T$ has at
least two orbits on each bipartite set of $\G$. It follows from~\cite[Lemma
3.2]{LWX} that $T$ is semiregular, a contradiction.

If $p=3$, then $q=2$ and $5^2\nmid|N|$ by Eq~(\ref{eq=2}),
implying that $N=T$ is simple.
Since $3^n\di |N|$ with $n\geq2$, we have
$N\cong A_6$ with $n=2$ or $\PSU(4,2)$ with $2\leq n\leq 4$.
For $n=2$, we have $|V(\G)|=2\cdot 2\cdot 3^2=36$
and $\G\cong\mathcal{G}_{36}$ by Example~\ref{exam=G36},
contrary to the hypothesis $\G\ncong\mathcal{G}_{36}$.
Suppose that $n=3$ or 4 and $N\cong\PSU(4,2)$.
Then $|V(\G)|=2^2\cdot 3^n$ and $|N|=2^6\cdot 3^4\cdot 5$.
Recall that $N$ has at most two orbits on $V(\G)$.
If $N$ is transitive on $V(\G)$, then
$|N_v|=2^{4}\cdot 3^{4-n}\cdot 5$
with $4-n=1$ or 0, and $N$ is arc-transitive on
$\G$.
It follows from Proposition~\ref{prop=stabilizer}
that $N_v\cong F_{20}\times\mz_4$.
However, $\PSU(4,2)$ has no subgroups isomorphic to
$F_{20}\times\mz_4$ by MAGMA~\cite{BCP}, a contradiction.
If $N$ has exactly two orbits on $V(\G)$,
then $\G$ is a bipartite graph with the
orbits of $N$ as its bipartite sets,
and $G$ has a $2$-element, say $g$, interchanging the two bipartite sets of $\G$.
Moreover, $|N_v|=2^5\cdot 3^{4-n}\cdot 5$
with $4-n=1$ or 0, and since $\PSU(4,2)$ has no subgroups of order $2^5\cdot 3\cdot 5$,
we have $|N_v|=2^5\cdot 5$.
Set $H=N\lg g\rg$. Then $H$ is arc-transitive on $\G$
and $|H_v|=2^{5+i}\cdot 5$ for some $i\geq0$,
which is impossible by Proposition~\ref{prop=stabilizer}.

\medskip

\f {\bf Case~2:} $|\pi(N)|=4$.

In this case, $T$ is a non-abelian simple
$K_4$-group, and $p\geq7$. Suppose that $q\neq5$. Then $5^2\nmid |N|$ by Eq~(\ref{eq=2}),
forcing that $N=T$ is a simple group. It follows from Lemma~\ref{prop=simple}~(2)
that $p^2\nmid|N|$, and since $p^n\di |N|$, we have $n=1$, which is contradict to the hypothesis $n\geq2$. Hence $q=5$, implying that $T$ is a $\{2,3,5,p\}$-group.
Again by Eq~(\ref{eq=2}),
we have $5^2\di |N|$, $5^3\nmid |N|$
and $3^3\nmid |N|$, forcing that $N=T$ or $N=T\times T$.

Assume that $N=T$. Since $5^2\di |N|$ and $p^2\di |N|$ with $p\geq7$, by Lemma~\ref{prop=simple}~(2) we have $N\cong \PSp(4,7)$ with $|N|=2^8\cdot 3^2\cdot 5^2\cdot 7^4$
and $(q,p,n)=(5,7,4)$,
or $\PSL(2,7^2)$ with $|N|=2^4\cdot 3\cdot 5^2\cdot 7^2$ and $(q,p,n)=(5,7,2)$.
Since $N$ has at most two orbits on $V(\G)$,
we have $|N_v|=2^7\cdot 3^2\cdot 5$ or $2^8\cdot 3^2\cdot 5$
for $N\cong \PSp(4,7)$, and $|N_v|=2^3\cdot 3\cdot 5$ or $2^4\cdot 3\cdot 5$
for $N\cong\PSL(2,7^2)$.
However, by MAGMA~\cite{BCP} neither $\PSp(4,7)$
nor $\PSL(2,7^2)$ has a subgroup of such order,
a contradiction.

Assume that $N=T\times T$. Since $3^3\nmid |N|$ and $5^3\nmid |N|$,
we have $N\cong \PSL(2,2^4)\times\PSL(2,2^4)$ with $(q,p,n)=(5,17,2)$,
$\PSL(2,11)\times \PSL(2,11)$ with $(q,p,n)=(5,11,2)$,
or $\PSL(2,31)\times \PSL(2,31)$ with $(q,p,n)=(5,31,2)$.
In this case, $N$ has a normal subgroup $T\cong\PSL(2,2^4)$,
$\PSL(2,11)$ or $\PSL(2,31)$. Since $|V(\G)|=2qp^2$,
$T$ is not semiregular on $V(\G)$, and since $p^2\nmid |T|$,
$T$ has more than two orbits on $V(\G)$.
Recall that $N$ has at most two orbits on $V(\G)$ and $5\di |N_v|$. If $N$ is transitive on $V(\G)$, then $N$ is arc-transitive
and Proposition~\ref{prop=atlesst3orbits} implies that $T$ is semiregular, a contradiction.
If $N$ has exactly two orbits on $V(\G)$, then $\G$ is bipartite with the orbits of $N$ as its bipartite sets. Clearly, $T$ has at
least two orbits on each bipartite set of $\G$,
and then~\cite[Lemma 3.2]{LWX} implies that $T$ is semiregular, a contradiction.

\medskip

\f {\bf Case~3:} $|\pi(N)|=5$.

In this case, $T$ is a non-abelian simple $K_5$-group, and $p>q\geq7$.
By Eq~(\ref{eq=2}), we have $5^2\nmid |N|$,
forcing that $N=T$ is simple.
It follows from Eq~(\ref{eq=2}) and Lemma~\ref{prop=simple}~(3) that
$p^2\nmid |N|$, and since $p^n\di |N|$,
we have $n=1$, a contradiction.

\medskip

Now, $N$ is solvable and thus it is an elementary abelian $r$-group for a prime $r$.
Clearly, $N$ has at least three orbits on $V(\G)$, and by Proposition~\ref{prop=atlesst3orbits}, $N$ is semiregular on $V(\G)$,
forcing that $|N|\di 2qp^n$.
Note that $p>q\geq 2$.
It follows that $N\cong \mz_2$, $\mz_2^2$ with $q=2$, $\mz_q$, or $\mz_p^s$ for some $s\leq n$.
To finish the proof, we may assume that $N\cong \mz_2$, $\mz_2^2$ with $q=2$,
or $\mz_q$ with $q\neq 2$, and aim to prove that $G$ has a non-trivial normal $p$-subgroup.
For $N\cong \mz_2$ with $q\neq 2$ (or $N\cong\mz_2^2$ with $q=2$, resp.), Proposition~\ref{prop=atlesst3orbits} implies that $\G_N$ is a connected graph of odd order $qp^n$ (or $p^n$, resp.) and odd valency 5, which is impossible. For $N\cong \mz_2$ with $q=2$
or $\mz_q$ with $q\neq 2$, $\G_N$ is a connected pentavalent $G/N$-arc-transitive graph of order $2p^{n}$. Let $B/N$ be a minimal normal subgroup of $G/N$. By Proposition~\ref{prop=minimal normal subgroup}, $B/N\cong \mz_p^m$,
and since $p>q$, $B$ has a normal Sylow $p$-subgroup $P$.
Since $B\unlhd G$ and $P$ is characteristic in $B$, we have $P\unlhd G$, as required.
\hfill\qed

Now, we are ready to prove Theorem~\ref{theo=2qp^n}.
\medskip

\f {\bf Proof of Theorem~\ref{theo=2qp^n}:}
Let $A=\Aut(\G)$, and let $M$ be a maximal normal $p$-subgroup of $A$
having more than two orbits on $V(\G)$. By Proposition~\ref{prop=atlesst3orbits}, $\G_M$ is a connected pentavalent $A/M$-arc-transitive graph of order $2qp^{n-s}$
for an integer $s\geq0$. Clearly, $n-s\geq0$.

If $n-s=0$, then $|V(\G_M)|=2q$.
It follows from Proposition~\ref{prop=2p}~(1) that
$\G_M\cong K_6$ with $q=3$, $K_{5,5}$ with $q=5$, or $\mathcal{CD}_q$ with
$5\di (q-1)$.

If $n-s=1$, then $|V(\G_M)|=2qp$, and $\G_M\cong K_{6,6}-6K_2$
with $(p,q)=(3,2)$, $\mathbf{I}_{12}$ with $(p,q)=(3,2)$, $\mathcal{G}_{42}$
with $(p,q)=(7,3)$,
$\mathcal{G}_{170}$ with $(p,q)=(17,5)$, $\mathcal{G}_{66}$
with $(p,q)=(11,3)$, $\mathcal{G}_{114}$ with $(p,q)=(19,3)$, $\mathcal{G}_{406}$
with $(p,q)=(29,7)$,
$\mathcal{G}_{3422}$ with $(p,q)=(59,29)$, $\mathcal{G}_{3782}$ with $(p,q)=(61,31)$, $\mathcal{G}_{574}$ with $(p,q)=(41,7)$,
or $\mathcal{CD}_{pq}$ with $5\di (p-1)$
and $q=5$ or $5\di (q-1)$ by Proposition~\ref{prop=2p}~(2).
Suppose that $\G_M\cong \mathcal{CD}_{pq}$.
By Example~\ref{exam=n}, $\Aut(\G_M)\cong D_{pq}\rtimes\mz_5$.
Since $A/M$ is arc-transitive on $V(\G_M)$,
we have $A/M=\Aut(\G_M)\cong D_{pq}\rtimes\mz_5$, implying
that $A/M$ has a normal Sylow $p$-subgroup of order $p$, say $B/M$.
Hence $B$ is a normal Sylow $p$-subgroup of $A$.
Since $|V(\G)|=2qp^n$ with $q<p$, $B$ has more than two orbits on $V(\G)$,
which is contradict to the maximality of $M$.
Hence $\G_M\ncong \mathcal{CD}_{pq}$.

If $n-s\geq2$, then Lemma~\ref{lem=4pn-minimal normal subgroup}  implies
that $\G_M\cong\mathcal{G}_{36}$ with $(p,q)=(3,2)$ or
$A/M$ has a minimal normal elementary abelian $p$-subgroup, say $B/M$.
For the former case, we are done. For the latter case,
the maximality of $M$ implies that $B/M$ has at most two orbits on $V(\G_M)$,
and thus $qp^{n-s}\di |B/M|$.
Since $B/M$ is a $p$-group, we have $p=q$, a contradiction.

Clearly, $K_6$, $K_{5,5}$ and $\mathcal{CD}_q$ are basic graphs because
each of them has order twice
a prime. Note that a pentavalent symmetric graph is basic if and only if its full automorphism group has no nontrivial normal subgroup having more than two orbits (see Proposition~\ref{prop=atlesst3orbits}).
For $K_{6,6}-6K_2$ and $\mathbf{I}_{12}$,
we have $\Aut(K_{6,6}-6K_2)\cong S_6\times \mz_2$ and $\Aut(\mathbf{I}_{12})\cong A_5\times \mz_2$ by Proposition~\ref{prop=2p}, both of which have a normal subgroup of order 2.
Hence neither $K_{6,6}-6K_2$ nor $\mathbf{I}_{12}$ is basic.
For the other graphs in Table~\ref{table=2},
Examples~\ref{exam=G42}-\ref{exam=G36}
and \ref{exam=G66} imply that
their full automorphism groups have no such nontrivial normal subgroup,
and hence they are basic graphs.\hfill\qed

\section{Pentavalent symmetric graphs of order $kp^n$ for some small $k$
and $n$}

In this section, we consider connected pentavalent symmetric graph of order $kp^n$ for some small
$k$ and $n$.

\subsection{$k=2$}

 First, we introduce some pentavalent symmetric graphs.

\begin{exam}\label{exam=G16}{\rm
Let $G=\lg a_1\rg\times \cdots \times \lg a_5\rg\cong\mz_2^5$,
and let $H=\lg a_1,a_2,a_3,a_4\rg$.
The Cayley graphs
$$Q_5=\Cay(G, \{a_1, \ldots, a_5\})~{\rm and}~FQ_4=\Cay(H,
\{a_1, a_2,a_3,a_4,a_1a_2a_3a_4\})$$
are called the 5-dimensional hypercube
and the 4-dimensional folded hypercube, respectively.
By~\cite[p.2649]{FZL},
we have $\Aut(Q_5)\cong \mz_2^5\rtimes S_5$
and $\Aut(FQ_4)\cong\mz_2^4\rtimes S_5$.}
\end{exam}

The following graph is a normal $\mz_2$-cover of $FQ_4$.

\begin{exam}\label{exam=G32}{\rm
Let $G_{32}=\lg a,b,c,d~|~a^4=b^2=c^2=d^2=[a,c]=[a,d]=1, a^b=a^{-1}, b^c=b^d=ba^2,
c^d=ca^2\rg$. Then $G_{32}\cong (D_8\rtimes\mz_2)\rtimes \mz_2$ is a non-abelian group of
order $32$. Define
$$\mathcal{G}_{32}=\Cay(G,\{b,ba,c,d,cda\}).$$
By MAGMA~\cite{BCP}, $\mathcal{G}_{32}$ is a connected pentavalent 2-transitive graph and $\Aut(\mathcal{G}_{32})\cong G_{32}\rtimes A_5$.}
\end{exam}

\begin{exam}{\rm \label{exam=64}
Let $G_{64}^{1}=\lg a,b,c,d,e~|~a^4=b^2=c^2=d^2=e^2=[a,c]=[a,d]=[a,e]=[b,e]=[c,e]=[d,e]=1, a^b=a^{-1}, b^c=b^d=ba^2,c^d=ca^2\rg$
and $G_{64}^{2}=\lg a,b,c,d,e~|~a^4=b^2=c^2=d^2=e^2=[a,b]=[a,c]=[a,d]=[a,e]=[b,c]=[c,e]=[d,e]=1, b^d=b^e=ba^2,c^d=ca^2\rg$
be two groups. Then $G_{64}^{1}\cong G_{32}\times\mz_2$ and
$G_{64}^{2}\cong (\mz_4\times\mz_2^2)\rtimes\mz_2^2$.
Moreover, both of them are non-abelian and have orders 64.
Define
$$\mathcal{G}_{64}^{1}=\Cay(G_{64}^{1},\{b,c,d,ab,acde\})~{\rm and}
~\mathcal{G}_{64}^{2}=\Cay(G_{64}^{2},\{b,c,d,e,abcde\}).
$$
By MAGMA~\cite{BCP}, $\mathcal{G}_{64}^{i}$ is a connected
pentavalent symmetric graph of order 64 for each $i=1$ or 2.
Moreover, $\Aut(\mathcal{G}_{64}^{1})\cong G_{64}^{1}\rtimes S_5$
and $\Aut(\mathcal{G}_{64}^{2})\cong G_{64}^{2}\rtimes D_5$.
}
\end{exam}

A list of all pentavalent $G$-arc-transitive graphs on up to 500 vertices with the vertex stabilizer $G_v\cong \mz_5$, $D_5$ or $F_{20}$ was given in magma code by Poto\v cnik~\cite{P}.
Based on this list, we have the following Proposition.

\begin{prop}\label{prop=500}
Let $\G$ be a connected pentavalent $G$-arc-transitive graph
of order $n$, where $G\leq \Aut(\G)$
and $n\in\{32,36,64,108,324\}$. If $G_v\cong \mz_5$, $D_5$
or $F_{20}$, then $\G\cong Q_5$, $\mathcal{G}_{32}$, $\mathcal{G}_{36}$,
$\mathcal{G}_{64}^{1}$, $\mathcal{G}_{64}^{2}$ or
$\mathcal{G}_{108}$.
\end{prop}

\begin{theorem}\label{cor=4p^2}
Let $p$ be a prime, and let $\G$ be a connected pentavalent symmetric graph of order $4p^n$
with $2\leq n\leq 4$. Then
$\G\cong FQ_4$, $\mathcal{G}_{36}$, $Q_5$, $\mathcal{G}_{32}$, $\mathcal{G}_{64}^1$,
$\mathcal{G}_{64}^2$, or $\mathcal{G}_{108}$.
\end{theorem}

\f {\bf Proof:} By Theorem~\ref{theo=2qp^n}, we have $p=2$ or 3.
Let $A=\Aut(\G)$, and let $M$ be a maximal normal $p$-subgroup of $A$
having more than two orbits on $V(\G)$. By Proposition~\ref{prop=atlesst3orbits}, $\G_M$ is a connected pentavalent $A/M$-arc-transitive graph of order $4p^{n-s}$
with $s\geq0$. Clearly, $n-s\geq1$.

Let $p=2$. Then $|V(\G)|=2^{n+2}$ with $4\leq n+2\leq 6$.
For $n+2=4$ or 5, we have $\G\cong FQ_4$, $Q_5$,
or $\mathcal{G}_{32}$ by~\cite[Proposition~2.9 and Theorem~1.1]{GHS}.
Assume that $n+2=6$. Then $|V(\G)|=64$.
By~\cite[Theorem~5.2]{FZL},
$A$ has a normal subgroup $N$ such that $\G_N\cong FQ_4$
or $Q_5$ and $A/N$ is arc-transitive
on $\G_N$. Since
each minimal arc-transitive subgroup of $\Aut(FQ_4)$ or
$\Aut(Q_5)$ has vertex stabilizer isomorphic to $\mz_5$ or $D_5$ by MAGMA~\cite{BCP},
$A/N$ has an arc-transitive subgroup $B/N$
such that $(B/N)_{\a}\cong\mz_5$ or $D_5$ with $\a\in V(\G_M)$. It follows that $B$ is arc-transitive on $\G$ and $B_v\cong (B/N)_{\a}\cong\mz_5$ or $D_5$ with $v\in V(\G)$.
By Proposition~\ref{prop=500}, $\G\cong \mathcal{G}_{64}^1$ or $\mathcal{G}_{64}^2$.

Let $p=3$. Then $|V(\G)|=4\cdot 3^n$ with $2\leq n\leq 4$.
By Theorem~\ref{theo=2qp^n}, $\G$ is a normal cover of $K_{6,6}-6K_2$, $\mathbf{I}_{12}$ or $\mathcal{G}_{36}$.
Hence $A$ has a normal subgroup $N$ such that $\G_N\cong K_{6,6}-6K_2$, $\mathbf{I}_{12}$, or $\mathcal{G}_{36}$. Moreover, $\G_N$ is $A/N$-arc-transitive.
By MAGMA~\cite{BCP}, under conjugate,
$\Aut(K_{6,6}-6K_2)$ has two minimal arc-transitive subgroups which
are isomorphic to $A_5\times \mz_2$ or $S_5$ respectively,
while both $\Aut(\mathbf{I}_{12})$ and $\Aut(\mathcal{G}_{36})$ have only one minimal arc-transitive subgroup which is isomorphic to $A_5$ and $A_6$, respectively.
Hence $A/N$ has an arc-transitive subgroup
$B/N$ such that $B/N\cong A_5\times \mz_2$ or $S_5$ with
$(B/N)_{\a}\cong D_5$ for
$K_{6,6}-6K_2$, or $B/N\cong A_5$ with $(B/N)_{\a}\cong\mz_5$ for $\mathbf{I}_{12}$,
or $B/N\cong A_6$ with $(B/N)_{\a}\cong D_5$ for $\mathcal{G}_{36}$,
where $\a\in V(\G_N)$.
It follows that $B_v\cong(B/N)_{\a}\cong\mz_5$ or $D_5$ for $v\in V(\G)$.
Since $B$ is arc-transitive on $\G$
and $|V(\G)|=36,$ $108$ or 324,
we have $\G\cong \mathcal{G}_{36}$, or $\mathcal{G}_{108}$ by Proposition~\ref{prop=500}.
\hfill\qed

\subsection{$k=3$}

Denoted by $\mathbf{I}_{12}^{(2)}$ the standard double cover of $\mathbf{I}_{12}$ (see~\cite[Theorem~4.1]{GZF}).
By MAGMA~\cite{BCP}, $\mathbf{I}_{12}^{(2)}$ is 2-transitive and
$\Aut(\mathbf{I}_{12}^{(2)})\cong A_5\rtimes D_4$.

\begin{theorem}\label{the=6p^2}
Let $p$ be a prime, and let $\G$ be a connected pentavalent symmetric graph of order $6p^2$. Then
$\G\cong \mathbf{I}_{12}^{(2)}$.
\end{theorem}

\f{\bf Proof:} Let $A=\Aut(\G)$. If $p=2$, then $|V(\G)|=24$, and by \cite[Theorem~4.1]{GZF},
$\G\cong \mathbf{I}_{12}^{(2)}$. Suppose that $p\geq 3$.
For $p=3$, \cite[Theorem~5.2]{FZL} implies that $\G$ is a normal cover of $K_6$,
and for $p>3$, Theorem~\ref{theo=2qp^n} implies that $\G$ is a normal cover of $K_6$
or $\mathcal{G}_{6p}$ with $p\in \{7,11,19\}$. To finish the proof,
we aim to find a contradiction.

Let $\G$ be a normal cover of $K_6$. Then $A$ has a normal subgroup
$N$ of order $p^2$ such that $\G_N\cong K_6$ and $\G_N$ is
$A/N$-arc-transitive. Since each arc-transitive subgroup of
$\Aut(K_6)$ has an arc-transitive subgroup isomorphic to $A_5$
by MAGMA~\cite{BCP}, $A/N$ has an arc-transitive subgroup $B/N\cong A_5$.
Moreover, $B$ is arc-transitive on $\G$, and both $B$ and $B'$ are insolvable.
It implies that $B'$ is not semiregular on $V(\G)$ because
otherwise $|B'|\di 6p^2$ and $B'$ is solvable. Hence $5\di
|B'_v|$ for $v\in V(\G)$, and by
Proposition~\ref{prop=atlesst3orbits}, $B'$ has at most two orbits
on $V(\G)$, forcing that $3p^2\di |v^{B'}|$ and $3\cdot 5\cdot
p^2\di |B'|$. Since $B/N\cong A_5$ and $\Mult(B/N)\cong \mz_2$,
Proposition~\ref{lem=m/n} implies that $p\di |B:B'|$, and thus $3\cdot 5\cdot
p^3\di |B|$. Since $|B|=60p^2$, we have $p=2$, a contradiction.

Let $\G$ be a normal cover of $\mathcal{G}_{6p}$ with $p=7$, 11 or 19.
Then $A$ has a normal subgroup $N$ such that $N\cong \mz_p$ and $\G_N\cong \mathcal{G}_{6p}$.
Furthermore, $\G_N$ is $A/N$-arc-transitive. By MAGMA~\cite{BCP}, each arc-transitive subgroup of $\Aut(\mathcal{G}_{6p})$
has only one minimal normal subgroup isomorphic to $\PSL(3,4)$, $\PSL(2,11)$ or $\PSL(2,19)$ for
$p=7$, $11$ or 19, respectively. Hence $A/N$ has a normal subgroup $B/N$ such
that $B/N\cong \PSL(3,4)$, $\PSL(2,11)$ or $\PSL(2,19)$.
Note that $\Mult(\PSL(3,4))\cong\mz_{4}\times\mz_{12}$
and $\Mult(\PSL(2,11))\cong\Mult(\PSL(2,19))\cong\mz_2$.
By Proposition~\ref{lem=m/n}, $B=B'N$ and $p\di|B:B'|$.
It implies that $p^2\nmid |B'|$ and $B'$ has at least three orbits on $V(\G)$.
Since $B'$ is characteristic in $B$ and $B\unlhd A$, we have $B'\unlhd A$, and thus
$B'$ is semiregular on $\G$ by Proposition~\ref{prop=atlesst3orbits}. It follows $|B'|\di 6p^2$, which is impossible because $B'\cong B/N$ is a non-abelian simple group.\hfill\qed

\medskip
\f {\bf Acknowledgement:} This work was supported by the National Natural Science Foundation of China (11571035), the 111 Project of China (B16002),
and the China Postdoctoral Science Foundation (2016M600841).

\end{document}